\documentclass[reqno,a4paper,twoside,11pt]{amsart}
\usepackage{amsthm}
\usepackage{graphicx}
\theoremstyle{definition}
\newtheorem{definition}{Definition}[section]

\theoremstyle{remark}

\theoremstyle{plain}

\newtheorem{theorem}[definition]{Theorem}

\newcommand{\init}{\mathrm{in}}
\newcommand{\set}[1]{\left\{{#1}\right\}}
\newcommand{\onemat}{\mathbb{I}}
\newcommand{\vek}[1]{\mathbf{#1}}

\begin{document}
\title{Bounds for the Entropy of Graded Algebras}
\author{Jan Snellman}
\address{Department of Mathematics\\
Stockholm University\\
SE-10691 Stockholm, Sweden}
\email{jans@matematik.su.se}
\keywords{Entropy of algebras, spectral radius, digraphs, \(0-1\)
  matrices, Perron-Frobenius theorem}
\subjclass{16W50, 05C50; 16P90, 05C20, 05C38}
\begin{abstract}
  Newman, Schneider and Shalev defined the \emph{entropy} of a graded
  associative algebra \(A\) as \(H(A) = \limsup_{n \to \infty}
  \sqrt[n]{a_n}\), where \(a_n\) is the vector space dimension of the
  \(n\)'th homogeneous component. When \(A\) is the homogeneous
  quotient of a finitely generated free associative algebra, they
  showed that \(H(A) \le \sqrt{a_2}\). Using some results of Friedland
  on the maximal spectral radius of \((0,1)\)-matrices with a
  prescribed number of ones, we improve on this bound.  
\end{abstract}
\maketitle
  Let \(A=\oplus_{n =0}^\infty A_n\) be a graded, associative algebra
  over a field \(F\). Assume furthermore that \(A\) is
  infinite-dimensional as a vector space 
  over \(F\), but that each graded component \(A_n\) is
  finite-dimensional, so that
   \(a_n = \dim A_n < \infty.\)
  Newman, Schneider and Shalev  \cite{EntropyGradedAlgebras} defined the
  \emph{entropy} of \(A\) by 
  \begin{equation}
    \label{eq:nsENT}
    H(A) = \limsup_{n \to \infty} \sqrt[n]{a_n}.
  \end{equation}
  Clearly, if \(a_n\) has polynomial growth, then \(H(A)=1\), and if
  \(a_n\) has exponential growth, with
  \(a_n \sim c d^n\), then \(H(A)=d\).
  Denote the Hilbert series of \(A\) by 
  \begin{equation}
    \label{eq:hilbert}
    A(t) = \sum_{n=0}^\infty a_n t^n.
  \end{equation}
  Then \(H(A)=1/R\), where \(R\) is the radius of convergence of
  \(A(t)\). If the Hilbert series is rational with
  \(A(t) = P(t)/Q(t)\), then \(\lvert\lambda_{\mathrm{min}}\rvert  \le
  H(A)  \le \lvert\lambda_\mathrm{max}\rvert\), where 
  \(\lambda_\mathrm{min}\) and \(\lambda_\mathrm{max}\)  are roots of
  \(Q(t)\) with minimal and maximal modulus. 

  We henceforth assume that \(A\) is the quotient \(A=F[X^*]/I\), where
  \(F[X^*]\) is the free associative algebra  on a finite alphabet
  \(X\), and \(I\) is a two-sided homogeneous ideal.  Then \(A\) is
  connected, i.e. \(A_0=F\), and generated by  
  \(A_1\). Newman, Schneider and Shalev  observed that the series
  \((a_n)_{n=0}^\infty\) is submultiplicative, i.e. \(a_n a_m \le
  a_{n+m}\), and that 
  this implies that 
  \begin{itemize}
  \item \(\lim_{n \to \infty} \sqrt[n]{a_n}\) exists and is bounded by
    \(\sqrt[n]{a_n}\) for each \(n \ge 1\), and consequently that
  \item \(H(A) \le \sqrt[n]{a_n}\) for all \(n \ge 1\).
  \end{itemize}

  In particular, 
  \begin{equation}
    \label{entSQ}
    H(A) \le \sqrt{a_2},
  \end{equation}
  and this bound is obtained
  when \(a_2=m^2\) and \(A\) is the free associative algebra on \(m\)
  letters. 

  Using a graph-theoretical result of Friedland \cite{Friedland:ME}
  and  Brualdi \& Hoffman \cite{BruHoff:SR}
  we can improve on this bound:

  \begin{theorem}\label{thm:entundir}
    Let \(a_2=k\), with
    \begin{displaymath}
      k = m^2 + \ell, \qquad 1 \le \ell \le 2m.
    \end{displaymath}
    Then 
    \begin{equation}
      \label{eq:udbound}
      H(A) \le \frac{m + \sqrt{m^2 + 2\ell}}{2} =: f(k)
    \end{equation}
    If \(\ell=1\) then 
    \begin{equation}
      \label{eq:udbound1}
      H(A) \le m,
    \end{equation}
    if \(\ell=2m -3\) then 
    \begin{equation}
      \label{eq:udboundm4}
      H(A) \le \frac{m -1 +  \sqrt{m^2 + 6m -7}}{2}.
    \end{equation}

    Let \(\rho_{m,p,q}\) denote the largest positive root of
    \begin{equation}
      \label{eq:rmqp}
      t^3 - m t^2 - pt + mp-pq
    \end{equation}
    Then, for any \(\ell\), there is an \(M(\ell)\) such that
    for \(m \ge M(\ell)\) we have that 
    \begin{equation}
      \label{eq:2}
      \dim A_2 = m^2 + \ell \quad \implies \quad 
      H(A) \le \rho_{m,\lfloor \ell/2 \rfloor, \lceil \ell/2 \rceil} 
    \end{equation}
  \end{theorem}
  \begin{proof}
    Recall that we've assumed that \(A=F[X^*]/I\). Denote by \(d\) the
    number of letters in \(X\), so that \(X=\set{x_1,x_2,\dots,x_d}\).
    Let \(\succ\) be a
    term-order on the free monoid \(X^*\), and let
    \(B=F[X^*]/\init(I)\), where \(\init(I)\) is the \emph{initial
      ideal} of \(I\) w.r.t \(\succ\). Since \(A(t)=B(t)\) we have
    that \(H(A)=H(B)\). Now, if we let \(J\) denote the monomial ideal
    generated by the quadratic part of \(\init(I)\), then
    \(C=F[X^*]/J\) is a quadratic monomial algebra. Clearly
    \(a_2=b_2=c_2\),  and \(A(t)=B(t)
    \ll C(t)\), i.e. \(C(t)=\sum_{n=0}^\infty c_n t^n\) and \(c_n \ge
    a_n=b_n\) for all \(n\). It follows that \(H(A)=H(B) \le H(C)\).

    Let \(G\) be the directed graph with vertex set
    \(\set{1,2,\dots,d}\), and with an (directed) edge from \(i\) to
    \(j\) iff 
    \(x_ix_j \not \in J\). Thus, we allow loops, but no multiple
    edges. 
    Clearly, for \(n \ge 2\) there is a bijection between on
    the one hand the 
    set of those 
    monomials in \(X^*\) that are  of length \(n\), and which do not
    belong to \(J\),  and on the other hand the set of directed walks
    in \(G\) of length \(n\). Thus, if \(g_n\) denotes the number of
    directed walks in \(G\) of length \(n\), then \(g_n = c_n\).

    Let \(M\) denote the adjacency matrix of \(G\), i.e. that \(d
    \times d\) matrix \((m_{ij})\) with \(m_{ij}=1\) iff there is an
    edge from \(i\) to \(j\). Let \(\vek{e} =(1,1,\dots,1) \) denote
    the (column) vector consisting of \(d\) ones. Then \(g_n =
    \vek{e}^t M^n \vek{e}\). If \[\phi(\lambda)=\lambda^d -h_1
    \lambda^{d-1} - h_2 \lambda^{d-2} - \cdots - h_{d-1}\lambda - h_d\] denotes the
    characteristic equation of \(M\), we have, by the Cayley-Hamilton
    theorem, that \(\phi(M)=0\). Hence 
    \begin{displaymath}
      M^d = h_1M^{d-1} + h_2 M^{d-2} + \cdots + h_{d-1}M + h_d I,
    \end{displaymath}
    which by multiplication with \(M^s\) yields
    \begin{displaymath}
      M^{s+d} = h_1M^{s+d-1} + h_2 M^{s+d-2} + \cdots + h_{d-1}M^{s+1} + h_d M^s,
    \end{displaymath}
    hence
    \begin{equation}\label{eq:rec}
      \begin{split}
      g_{s+d} &= \vek{e}^t M^{s+d} \vek{e} \\
      &= h_1 \vek{e}^t M^{s+d-1} \vek{e}+
      h_2 \vek{e}^t M^{s+d-2} \vek{e}+ \cdots + h_{d-1}
      \vek{e}^tM^{s+1} \vek{e}+ h_d \vek{e}^tM^s\vek{e} \\
      &=
      h_1 g_{s+d-1} + h_2 g_{s+d-2} + \cdots + h_{d-1} g_{s+1} + h_d
      g_s
      \end{split}
    \end{equation}
    Since the \(g_n\)'s and hence the \(c_n\)'s satisfy the linear
    recurrence \eqref{eq:rec}, it follows (as
    mentioned in section 2 of \cite{EntropyGradedAlgebras}) that 
    \(H(C) \le \left \vert \rho \right \rvert\),
    where \(\rho\) is the  maximal modulus of a root of
    the characteristic polynomial of \(M\), i.e. the largest absolute
    value of an eigenvalue
    of \(M\), i.e. the spectral radius
    of \(M\) (and by definition of \(G\)).

    It follows that we can  apply the bounds obtained by Friedland
    \cite{Friedland:ME} and Brualdi \& Hoffman \cite{BruHoff:SR} for
    the spectral radius of directed graphs with \(k\) edges to bound
    the entropy of \(C\), and hence that of \(A\).

  Let \(G_{m,p,q}\) be the directed graph on \(\set{1,\dots,m+1}\),
  where there is an edge from \(i\) to \(j\) if \(i,j \le m\) or if \(i
  \le p\) and \(j = m+1\) or if \(i=m+1\) and \(j \le q\).
  Thus, \(G_{m,p,q}\) has adjacency matrix 
  \begin{displaymath}
    \left[
    \begin{array}{cccccc|c}
     &     & &&& &1 \\
     &     & &&& &1 \\
     & & \onemat_m   &&& &\vdots \\
     &     & &&& &1 \\
     &     & &&& &0 \\
     &     & &&& &\vdots \\
     \hline
   1 & 1 & \cdots & 1 & 0 & \cdots & 0
    \end{array}
    \right]
  \end{displaymath}
  with \(p\) ones on the last row, \(q\) ones on the last column, and
  with \(\onemat_m\) an \(m \times m\) matrix of ones.
  It corresponds to the quadratic monomial algebra 
  \begin{displaymath}
    A_{m,p,q} =\left
    \langle x_1,\dots,x_m \, \rvert \, x_{m+1}^2; x_m x_{m+1},\dots,
    x_{p+1} x_{m+1};  x_{m+1}x_m,\dots x_{m+1}x_{q+1} \right \rangle.
  \end{displaymath}

  Let \(\rho_{m,p,q}\) denote the spectral radius of \(G_{m,p,q}\),
  hence \(\rho_{m,p,q}=H(A_{m,p,q})\). 
  Friedland \cite{Friedland:ME} showed that 
  \(\rho_{m,p,q}\) is the largest positive root of \eqref{eq:rmqp}
  and that for \(k=m^2+1\), the graph \(G_{m,0,1}\) has maximal
  spectral radius \(\rho_{m,0,1}=m\) of all directed graphs with \(k\)
  edges (the later result was also obtained by Brualdi and Hoffman
  \cite{BruHoff:SR}). This gives \eqref{eq:udbound1}. Friedland also
  showed that for 
  \(k=m^2+2m-3=(m+1)^2-4\) edges, the maximal spectral radius is obtained
  \emph{not} by \(G_{m,m-2,m-1}\) but by the graph on \(m+1\) vertices
  where there is an edge from \(i\) to \(j\) if either \(i<m\) or
  \(j<m\); i.e. by the digraph with adjacency matrix
  \begin{displaymath}
    \left[
    \begin{array}{cccc|cc}
     &      &&&1 &1 \\
     &      &&&1 &1 \\
     & & \onemat_{m-1}   && \vdots &\vdots \\
     &     & &&1 &1 \\
     \hline
   1 & 1 & \cdots & 1 & 0 &  0 \\
   1 & 1 & \cdots & 1 & 0 &  0 \\
    \end{array}
    \right]
  \end{displaymath}
  Friedland showed that this graph has spectral radius
  \[\frac{m-1+\sqrt{m^2+6m-7}}{2},\] so we have \eqref{eq:udboundm4}.

  Friedland furthermore showed that for a fixed \(\ell\) it holds that for
  sufficiently large \(m\), if \(k=m^2+\ell\) then 
  \(G_{m,\lfloor  \ell/2 \rfloor,\lceil \ell/2 \rceil}\) has the
  largest spectral radius
  of a directed graph with \(k\) edges.
  This gives \eqref{eq:2}.
  \end{proof}

  We have that for \(k=m^2 + \ell\), 
  \begin{multline*}
f(k) < \sqrt{k} \quad \iff \quad
  f(k)^2 < k  \quad \iff \quad \left(\frac{m + \sqrt{m^2 +
            2\ell}}{2} \right)^2 < m^2+\ell \\  \iff \quad
m^2 + \ell - m \sqrt{m^2+2\ell} >0 \quad \iff \quad  \ell >0,    
  \end{multline*}
 so the bound of \eqref{eq:udbound} is indeed an
 improvement for \(\ell > 0\). The number \(\rho_{m,\lfloor \ell/2
   \rfloor, \lceil \ell/2 \rceil}\) is smaller still, as is shown in
 the picture below, which plots \(\sqrt{k}-f(k)\) (in a solid line) and
 \(\sqrt{k}-\rho_{m,\lfloor \ell/2 \rfloor, \lceil \ell/2 \rceil}\)
 (in a dotted line)
 for \(5 \le k \le 85\).  
 
 
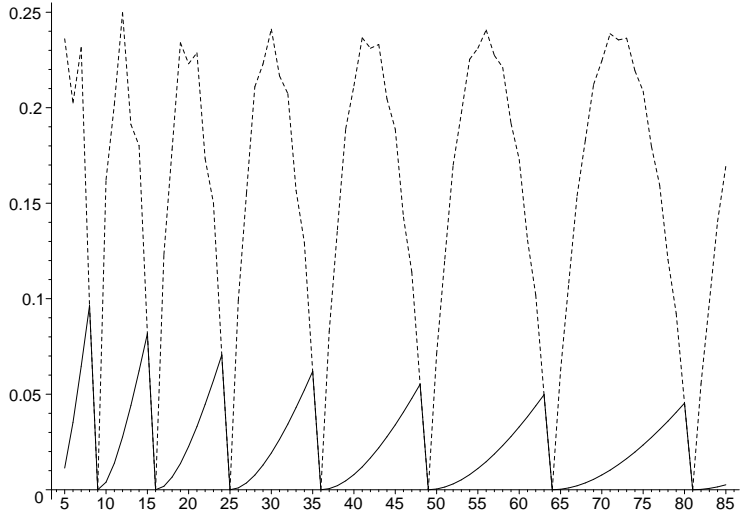
\begin{center}
   \includegraphics*[scale=0.5, angle=270, bb=100 100 510 670 ]{enplot1.ps}
 \end{center}


\bibliographystyle{plain}
\bibliography{journals,entropy}

\begin{thebibliography}{1}

\bibitem{BruHoff:SR}
R.~A. Brualdi and A.~J. Hoffman.
\newblock On the {S}pectral {R}adius of {(0,1)}-{M}atrices.
\newblock {\em Linear {A}lgebra and its {A}pplications}, 65:133--146, 1985.

\bibitem{Friedland:ME}
Schmuel Friedland.
\newblock The {M}aximal {E}igenvalue of {0-1} {M}atrices with {P}rescribed
  {N}umber of {O}nes.
\newblock {\em Linear {A}lgebra and its {A}pplications}, 69:33--69, 1985.

\bibitem{EntropyGradedAlgebras}
M.~F. Newman, Csaba Schneider, and Aner Shalev.
\newblock The {E}ntropy of {G}raded {A}lgebras.
\newblock {\em Journal of {A}lgebra}, 223:85--100, 2000.

\end{thebibliography}
\end{document}